\documentclass[11pt]{article} 

\usepackage[utf8]{inputenc} 
\usepackage{amssymb,amsmath}
\usepackage{amsfonts}

\usepackage{geometry} 
\usepackage{graphicx} 

\usepackage{booktabs} 
\usepackage{array} 
\usepackage{paralist} 
\usepackage{verbatim} 
\usepackage{subfig} 
\usepackage{ragged2e}
\usepackage{enumerate}
\usepackage{fancyhdr} 
\usepackage{sectsty}
\allsectionsfont{\sffamily\mdseries\upshape} 

\usepackage[nottoc,notlof,notlot]{tocbibind} 
\usepackage[titles,subfigure]{tocloft} 

\title{Finite Commutative Semihypergroups Built From Groups}
\author{Stan Onypchuk}

\begin{document}
\maketitle
\begin{abstract}
\begin{justify}
Necessary and sufficient conditions for finite commutative semihypergroups to be built from abelian groups of the same order are established.
\end{justify}
\end{abstract}

\noindent
\textbf{Introduction} The theory of hypergroups and semihypergroups was introduced by C .Dunkl [1], R. Jewett [2] and R. Spector [3] and is well developed now.\\
Many examples of finite commutative semihypergroups and hypergroups we can find in [4]. In [4] there is a precise \textsl{physical} definition of a finite semihypergroup:\\
\begin {justify}
\textsl{A finite commutative hypergroup is a finite collection of particles, say \{$c_{0}, \ldots , c_{n}$\}, which are allowed to interact by colliding. When two particles collide, they coalesce to form a third particle. The results of collisions are however not definite; if we collide $c_{i}$ with $c_{j}$ the probability of emerging with the single particle $c_{k}$ is $n_{ij}^{k}$ and is fixed.
}
\end {justify}
In the present article we study some finite commutative semihypergroups which can be developed from commutative groups. In terms of the above physical interpretation the main result shows that for such semihypergroups each particle in \{$c_{0}, \ldots , c_{n}$\} is not an elementary particle, but is a combination of some different elementary particles \{$b_{0}, \ldots , b_{n}$\} and evolution of \{$b_{0}, \ldots , b_{n}$\} (collision of $b_{i}$ and $b_{j}$) can be described as a group operation for some group. Each such semihypergroup defines underlying group precisely.\\
\\
Let $\mathbf{H}$ be a finite semihypergroup  with $n$ states $e_{1},e_{2}, \ldots , e_{n} $ and convolution operation defined by 
\begin{equation}
			e_{i} \ast e_{j} = \sum_{k} a_{i,j}(k)e_{k}  \text{  }    (i,j,k = 1,2, \ldots , n) 
\end{equation}			                      
where $a_{i,j}(k)\geqslant 0 \text{ and } \sum_{k = 1}^{n} a_{i,j}(k) = 1$  for each $i,j$ is associative 
\[
	(e_{i} \ast e_{j}) \ast e_{k} = e_{i} \ast (e_{j} \ast e_{k}) = e_{i} \ast e_{j} \ast e_{k} 
\]
Let us denote \\
columns \{$ a_{i,j}(1), \ldots , a_{i,j}(n)$\}  by $a_{i,j}$ \\
matrix with columns \{$a_{i,1}, \ldots , a_{i,n}$\}  by $A_{i} \text{ (matrix of left regular representation)} $\\
matrix with columns \{$a_{1,i}, \ldots , a_{n,i}$\}  by $B_{i} \text{ (matrix of right regular representation)}$ \\
and cube with matrices \{$A_{1}, \ldots , A_{n}$\} by $C.$ \\
\\
\textbf{Definition} A commutative semihypergroup $\mathbf{H}$ will be said to be \emph{derived from a group} if the following condition is satisfied\\
(A) There are only $n$ different columns among all \{$a_{i,j}$\} and columns in each $A_{i}$ and $B_{i}$ are linearly independent.\\
\\
In matrix terms the convolution operation (1) transforms to $ a_{i,j} = A_{i}e_{j}$ (where $ e_{j} $ is a column-vector with 1 on $j$-th position and 0 otherwise). The matrix $\sum c_{i}A_{i} $ represents a measure $\sum c_{i}e_{i} $ so $\sum_{k=1}^{n} a_{i,j}(k)A_{k}$ represents $e_{i} \ast e_{j}$. Now it is obvious that the convolution operation (1) is associative if and only if
\begin{equation}
	A_{i}A_{j} = \sum_{k=1}^{n} a_{i,j}(k)A_{k} \\
\end{equation}
Now let us establish some feature all derived from a group semihypergroups are hold:\\
\textbf{Corollary 1.} There are no more than $n$ different elements in $C$ and all columns in $C$ contain the same elements.\\
The $j$-th column in $A_{k}A_{k}$ is \\
\[
	\sum a_{i,j} a_{k,k}(i) 
\]
So
\begin{equation}
	\sum a_{k,i}a_{k,j}(i) = \sum a_{i,j}a_{k,k}(i)
\end{equation}
Because \{$a_{k,i}$\} ($k$ - fixed) and  \{$a_{i,j}$\} ($j$ - fixed) are from the same set of $n$ linearly independent columns in $C$ we have that columns $a_{k,j}$ contain the same elements as the column $a_{k,k}$. Now the statement is following from commutativity of $\mathbf{H}$.\\

From (3) when $i = j$ in the left part and $i = k$ in the right part we have \\
\[
	a_{k,j}(j) = a_{k,k}(k) 
\]
It leads to\\
\textbf{Corollary 2.} Diagonal elements in all matrices $A_{i}$ are equal.\\
\textbf{Corollary 3.} All rows in all matrices $A_{i}$ contain the same elements as the columns.\\
Let us assume that all not zero elements in $C$ are different and $a$ is one of them. By  \textsl{Corollary 1} there are exactly $n$ elements $a$ in $A_{i}$. If one of the rows in $A_{i}$ does not contain element $a$ then there exists the row $a_{i,j}(k)$ ($j = 1,\cdots ,n$) with two such elements. Let $a_{i,j_{1}}(k) = a_{i,j_{2}}(k) = a$. Then in both columns $A_{i}A_{i}(j_{1})$ ($j_{1}$-th column of $A_{i}A_{i}$) and $A_{i}A_{i}(j_{2})$ ($j_{2}$-th column of $A_{i}A_{i}$) the column $a_{i,k}$ will be added with the same coefficient $a$. But by (2) $A_{i}A_{i}(j_{1}) = \sum_{m} a_{i,i}(m)a_{m,j_{1}}$ and $A_{i}A_{i}(j_{2}) = \sum_{m} a_{i,i}(m)a_{m,j_{2}}$ and by \textsl{Corollary 1} there exists $m_{1}$ such that $a_{i,i}(m_{1}) = a$. Therefore $a_{m_{1},j_{1}} = a_{m_{1},j_{2}} = a_{i,k}$. But this contradict to (A). 
\\
\textbf{Conclusion}\\
there are only $n$ different elements in the cube $C$,\\
there are only $n$ different columns and rows in the cube $C$,\\
each column and row in the cube $C$ contains the same set of elements.\\
\\
Using these features of derived from group semihypergroups we can show\\
\textbf{Corollary 4.} Every matrix $A_{i}$ is a linear combination of $G_{i}$, where \{$G_{i}$\} - matrices of a regular representation of some abelian group $g$ of order $n$.\\
To prove this we will use the following alternative definition of a group:\\
\\
\hspace*{0.3in}\textbf{a)} On non-empty set $G$ of elements \{$g_{i}$\} there is defined a binary operation - $g_{i} \cdot g_{j} = g_{k}$;\\
\hspace*{0.3in}\textbf{b)} This operation is associative;\\
\hspace*{0.3in}\textbf{c)} For any elements $g_{i}$ and $g_{j}$ there exist at least one such $g_{k}$ and at least one such $g_{l}$ that
\centerline{	$g_{k} \cdot g_{i} = g_{i} \cdot g_{l} = g_{j}.$}
\\
Let us assume that all $n$ elements are different and one of them is $a$.\\
Then for any $i,k$ there exists such $j$ that
\begin{equation}
	a_{i,j}(k) = a
\end{equation}
and for any $j,k$ there exists such $i$ that
\begin{equation}
	a_{i,j}(k) = a
\end{equation}
When $a$ = 1 and all other elements in the cube $C$ are equal to 0, conditions (4) and (5) describe binary operation defined as
\begin{equation}
	g_{i} \cdot g_{j} = a_{i,j}(k) = g_{k}
\end{equation}
Let's show that the operation (6) is associative.
\begin{equation}
	(e_{i} \ast e_{j}) \ast e_{m} = \sum_{k=1}^{n} a_{i,j}(k)(e_{k} \ast e_{m}) = \sum_{k=1}^{n}a_{i,j}(k)a_{k,m}
\end{equation}
on the other hand
\begin{equation}
	e_{i} \ast (e_{j} \ast e_{m}) = \sum_{p=1}^{n} a_{j,m}(p)(e_{i} \ast e_{p}) = \sum_{p=1}^{n}a_{j,m}(p)a_{i,p}
\end{equation}
Because columns in (7) and (8) are independent and convolution operation is associative, we have to conclude that $a_{k,m} = a_{i,p}$ when $a_{i,j}(k) = a_{j,m}(p) = 1$ . That means that
\[
	(g_{i} \cdot g_{j}) \cdot g_{m} = g_{i} \cdot (g_{j} \cdot g_{m})
\]
So, for any element $a$ in this cube $C$ limitation to construction when it's 1 in such places and is 0 otherwise is a group.\\
When $n$ is a prime number (and some other cases) there exists only one group of order $n$ - cyclic group $\mathbb{Z}_{n}$. It is well known that an abelian group of order n is a direct product of cyclic groups of prime order.\\
Suppose when cube's elements equal $a$, they correspond to group $G_{1}$ and when they equal $b$, then they correspond to group $G_{2}$. All groups for our $C$ are among direct product of $\mathbb{Z}_{p} \times D$ where $p$ is prime and $D$ has order $m \text{ such that } pm = n$.\\
Consider the $p \times p$ matrix
\[
	\mathbf{Z}_{2} =
	\begin{pmatrix}
		0 & 0 & \ldots & 1 \\
		1 & 0 & \ldots & 0 \\
		0 & 1 & \ldots & 0 \\
		\ldots & \ldots & \ldots \\
		\ldots & \ldots & \ldots \\
		\ldots & \ldots & \ldots \\
		0 & 0 & \ldots & 0
	\end{pmatrix}
\]
and build the matrix $G_{22}$ by replacing in $\mathbf{Z}_{2}$ every element 0 by $m \times m$ zero matrix and every element 1 by $m \times m$ identity matrix. Obviously such formed matrix  $G_{22}$ belongs to representation of $\mathbb{Z}_{n}$ and representation of any abelian group from $\mathbb{Z}_{p} \times D$.\\
Assume that there are two numbers $a \text{ and } b$ and two different groups with its matrix representations $G_{1} \text{ and } G_{2}$. Consider the following two fusions of their matrices \\
\begin{equation}
	G_{1,i_{1}}, \ldots , G_{1,i_{k}}, \ldots , G_{1,i_{n}} 
\end{equation}
\begin{equation}
	G_{1,j_{1}}, \ldots , G_{22}, \ldots , G_{1,j_{n}} 
\end{equation}
\begin{equation}
	G_{2,j_{1}}, \ldots , G_{22}, \ldots , G_{2,j_{n}}
\end{equation}

Here in (9) and (10) we have matrices from representation of the same group in different order when in pair (10) and (11) there are matrices from representation of different groups but in the same order.
Now, when we add sets (9) and (10) by component all matrices will have the same set of columns, particularly, the sum of $G_{1,i_{k}} \text{ and } G_{22}$. But the sum of $G_{1,i_{k}} \text{ and } G_{22}$ (now $G_{22}$ is from (11)) has the same set of columns. The same $n$ different columns that any matrix from sum (9) and (10) has. Now consider a different position $j_{p}$ where matrices $G_{1,j_{p}} \text{ and } G_{2,j_{p}}$ are different. Then  $G_{1,i_{P}} + G_{1,j_{p}}$ and  $G_{1,i_{p}} + G_{2,j_{p}}$ will have different sets of columns. It means that sum of (9) and (11) has more then $n$ different columns which contradicts condition (A). \\
\\
\textbf{Main result.}\\
Commutative semihypergroup $\mathbf{H}$ with convolution  $e_{i} \ast e_{j}$ is derived from group if and only if there exists abelian group $g$ of order $n$ and some probability measure $m$ on $g$ that any $e_{i}$ is\\
\[
e_{i} = m \cdot g_{i}
\]
where $g_{i} \cdot g_{j}$ is a group product on $g$. \\
Let us show that every element in $\mathbf{H}$ is a combination of elements of an abelian group $g$ of order $n$. \\
Let $G$ be a matrix representation of $g$ and $G_{i}$ be the $i$-th matrix in $G$. Assume that 
\begin{equation}
	A_{1} = \sum_{k=1}^{n}m_{k}G_{k}   \text{       where  (}\sum m_{k} = 1\text{)}
\end{equation}
but $\mathbf{H}$ is commutative hence $A_{1} = B_{1} \text{ and } b_{1,i} = a_{i,1} = a_{1,i}$.  By (12)
\[
	a_{1,i} = \sum m_{k}g_{k,i}
\]

where $g_{k,i}$ is $i$-th column in $G_{k}$. But $b_{1,i}$ represents the $i$-th element in $\mathbf{H}$ and $g_{k,i}$ represents the ($g_{k}   \cdot g_{i}$)-th element in in $g$. Then
\[
	e_{i} = \sum m_{k}\text{(}g_{k} \cdot g_{i} \text{)} = ( \sum m_{k}g_{k}) \cdot g_{i} = m \cdot g_{i}
\]
Now let us show that system \{$e_{1}, \ldots , e_{n}$\} where $e_{i} = m \cdot g_{i}$  with convolution operation
\[
e_{i} \ast e_{j} = m \cdot g_{i} \cdot m \cdot g_{j} = m \cdot (g_{i} \cdot m \cdot g_{j}) = m \cdot (\sum m_{i,j}(k) g_{k}) = \sum  m_{i,j}(k) m \cdot g_{k} = \sum m_{i,j}(k) e_{k}
\]
describes some semihypergroup $\mathbf{H}$ which satisfies condition (A).\\
Let consider how the $i$-th matrix $M_{i}$ of a regular representation of $\mathbf{H}$ is built.\\
Each $m_{i,j}$ the $j$-th column in $M_{i}$ is the ordered sequence of coefficients in $g_{i}\cdot m \cdot g_{j}$ - \{$m_{i,j}(k)$\}.\\
Because 
\[ 
	m \cdot g_{j} = m_{1}g_{1} \cdot g_{j} + m_{2}g_{2} \cdot g_{j} + \ldots + m_{n}g_{n} \cdot g_{j}
\]
and $g_{i} \cdot g_{j}$ is the $j$-th column in $G_{i}$ we have to conclude that $m \cdot g_{j}$ represents the column $\sum m_{k} g_{k,j}$ and $m\cdot g_{j}  (j = 1, \ldots , n)$ should be represented by matrix 
\[
\sum m_{k} G_{k} = M = 
\begin{pmatrix}
		m_{1,1}(1) & m_{1,2}(1) & \ldots & m_{1,n}(1) \\
		m_{1,1}(2) & m_{1,2}(2) & \ldots & m_{1,n}(2) \\
		\ldots & \ldots & \ldots \\
		\ldots & \ldots & \ldots \\
		\ldots & \ldots & \ldots \\
		m_{1,1}(n) & m_{1,2}(n) & \ldots & m_{1,n}(n)
	\end{pmatrix}
\]
where $m_{1,1}(1) = m_{1}, m_{1,1}(2) = m_{2}, \ldots ,m_{1,1}(n) = m_{n}$.
Now when we add scalar multiplier $g_{i}$ to $M$ we need to implement the multiplication by component $g_{i}$ on $M$. It means that each $k$-th row in $M$ ($m_{1,1}(k), m_{1,2}(k), \ldots ,m_{1,n}(k)$) will be \{$m_{1,1}(k)g_{i}\cdot g_{k}, m_{1,2}(k)g_{i}\cdot g_{k}, \ldots ,m_{1,n}(k)g_{i}\cdot g_{k}$\} will be moved to $p$-th row where $g_{p} = g_{i} \cdot g_{k}$.\\
It exactly means that $M_{i} = G_{i} M$ and because group g is commutative $M_{i} = M G_{i}$. In $M G_{i}$  $ G_{i}$ acts as permutation of columns in $M$ that means that any $M_{i}$ has the same set of columns as $M$.\\  
Commutativity here is essential. Next, if measure $m$ is not an uniform measure concentrated on some subgroup of g all columns in $M$ will be different and linearly independent.\\
Semihypergroup $H$ is associative just because underlying group $g$ is associative and $e_{i} \ast e_{j} = (m \cdot g_{i}) \cdot (m \cdot g_{j}) = m \cdot (g_{i} \cdot m \cdot g_{j})$. \\ 
\\
\textbf{References}\\
\text{[1]} C. F. Dunkl, The measure algebra of a locally compact hypergroup, Trans. Amer. Math. Soc. 179(1973) 331-348.\\
\text{[2]} R. I. Jewett, Spaces with an abstract convolution of measures, Advances in Math. 18 (1975) 1-110.\\
\text{[3]} R. Spector, Mesures invariantes sur les hypergroupes,Trans Amer Math Soc 239, (1978) 147 -n 165.\\
\text{[4]} N. J. Wildberger, Finite commutative hypergroups and applications from group theory to conformal field theory, Contemporary Mathematics, 188 (1995) 413-431. 

\end{document}